\documentclass[12pt]{article}
\usepackage{amsthm}
\usepackage{amsmath}
\usepackage{amssymb}
\usepackage{amsfonts}
\usepackage{latexsym}

\newcommand{\nc}{\newcommand*}
\newcommand{\rnc}{\renewcommand*}

\nc{\ts}{\textstyle}

\def\half{{\frac{1}{2}}}

\nc{\ket}[1]{{\vert{#1}\rangle}}                
\nc{\bra}[1]{{\langle{#1}\vert}}                
\nc{\braket}[2]{{\langle{#1}\vert{#2}\rangle}}  
\nc{\scpr}[2]{\left({#1}\,,\,{#2}\right)}      
\nc{\ccr}[2]{\left[{#1}\,,\,{#2}\right]}       
\nc{\car}[2]{\left\{{#1}\,,\,{#2}\right\}}     
\nc{\pois}[2]{\left\{{#1},{#2} \right\}_{cl}}  
\nc{\Szi}[1]{\sum_{ #1 = 0 }^{\infty}}
\nc{\Szx}[2]{\sum_{#1 =0}^{#2}}
\nc{\Syx}[3]{\sum_{ #1 = #2 }^{#3}}
\nc{\Izx}[3]{\int_{ 0 }^{#1}{#2}{\rm d}\,{#3}}
\nc{\qIzx}[3]{\int_{ 0 }^{#1}{#2}{\rm d}_q\,{#3}}
\def\qd{\left(\ts\frac{\text{d}}{\text{d}x}\right)_q\,}
\nc{\goth}{\mathfrak}
\rnc{\bold}{\mathbf}
\rnc{\frak}{\mathfrak}
\rnc{\Bbb}{\mathbb}
\rnc{\rm}{\text}
\nc{\mbf}[1]{\mbox{\boldmath ${#1}$}}
\nc{\BR}{\mathbb R}
\nc{\BC}{\mathbb C}
\nc{\BI}{{\mathbb I}}

\nc{\gmf}[1]{\varGamma({ #1})}
\nc{\Foo}[3]{\raisebox{-3pt}{${}_1\text{\large F}_1$}%
 \left(\ts{\genfrac{}{}{0pt}{}{#1}{#2}\left.\right| #3}\right)}
\nc{\Fto}[4]{\raisebox{-3pt}{${}_2\text{\large F}_1$}%
 \left(\ts{\genfrac{}{}{0pt}{}{#1, #2}{#3}\left.\right| #4}\right)}
\nc{\Ftf}[7]{\raisebox{-3pt}{${}_2\text{\large F}_4$}%
 \left(\ts{\genfrac{}{}{0pt}{}{#1, #2}{#3, #4, #5, #6}%
 \left.\right| #7}\right)}
\nc{\sfito}[4]{{{}_2{\phi}_1}%
 \left(\ts{\genfrac{}{}{0pt}{}{#1,\,\, #2}{\phantom{\bigr(}#3}%
 \biggl.\biggr| #4 }\right)}
\nc{\sfioo}[3]{{{}_1{\phi}_1}%
 \left(\ts{\genfrac{}{}{0pt}{}{#1\phantom{\bigr(}}{#2\phantom{\bigr(}}%
 \!\!\biggl.\biggr| #3 }\right)}
\begin{document}
$${}$$
\begin{flushright}
{submited to the Proceedings}\\ {of the international seminar}\\ {{\it
DAY} on {\it DIFFRACTION}' 2003}
\end{flushright}
$${}$$

\centerline{\huge\bf Generalized coherent states}

\medskip

\centerline{\large\bf for q-oscillator connected with}
\bigskip

\centerline{{\huge\bf q-Hermite Polynomials} \large\footnote{This
research is supported in parts by RFFI grants No 03-01-00837, 03-01-00593.}}

\bigskip
\smallskip
{\flushleft{\large\bf V.V.Borzov,${}^*$
E.V.Damaskinsky${}^{**}$}}
\smallskip
{\flushleft{${}^*$ St.-Petersburg University of
Telecommunications\\ E-mail: vadim@VB6384.spb.edu}} \vspace{-.1cm}
{\flushleft{${}^{**}$St.-Petersburg University of the Defence
Engineering Constructions\\ E-mail: evd@pdmi.ras.ru}}
\bigskip

\centerline{\bf Abstract}
\begin{quote}
For the oscillator-like systems, connected with q-Hermite
polynomials, coherent states of Barut-Girardello type are defined.
The well-known Arik-Coon oscillator naturally arose in the
framework of suggested approach as oscillator, connected with the
Rogers q-Hermite polynomials, in the same way as usual oscillator
with standard Hermite polynomials. The results about the coherent
states for discrete q-Hermite polynomials of II type are quite
new.
\end{quote}
\vskip 0.5cm

\section{Introduction}
We consider the oscillator-like systems (or generalized
oscillators), connected with q-Hermite polynomials in the same way
as usual boson oscillator connected with the standard Hermite
polynomials. We define the analogues of coherent states of
Barut-Girardello type for these generalized oscillators. One year
ago we talked here about GCS connected with classical polynomials.
Our approach to such construction is developed in the following
papers~\cite{bdk}~-~\cite{bd6}.

Within of Askey - Wilson scheme~\cite{AW},~\cite{KS}, we know
$q$-Hermite polynomials of three kind only
\begin{enumerate}
\item The {\bf Rogers $\mathbf q$-Hermite polynomials }
\begin{equation}
\label{qRH} H_n(x;q)=\sum_{k=1}^{n}
\frac{(q;q)_n}{(q;q)_k(q;q)_{n-k}}
\exp\left\{i(n-2k)\theta\right\},
\end{equation}
where $x=\cos\theta, \, 0<q<1$ and
\begin{equation}
\label{qfac}
(q;q)_k=\prod_{s=1}^{k}(1-aq^{s-1}), \quad
(q;q)_{\infty}=\prod_{s=1}^{\infty}(1-aq^{s-1}).
\end{equation}
\item The {\bf discrete $\mathbf q$-Hermite polynomials of I-type}
\begin{equation}\label{qHI}
h_n(x;q)=q^{\binom{n}{q}} \sfito{q^{-n}}{x^{-1}}{0}{q;\,-qx}
\end{equation}
\item The {\bf discrete $\mathbf q$-Hermite polynomials of II-type}
\begin{equation}\label{qHII}
{\tilde h}_n(x;q)=x^n \sfito{q^{-n}}{q^{-n+1}}{0}{q^2;\,-x^{-2}}.
\end{equation}
\end{enumerate}

In the first case we get the well-known q-oscillator of
Arik - Coon~\cite{AC}. We want to stress that in this case we obtain these
well-known results using a new approach.

In the second case, it turned out, that no GCS exist. This is due
to the fact, that the radius of convergence
for normalizing factor of GCS equals zero.

Finally, in the third case we will construct GCS. These results are
quite new. Unfortunately, we do not know the measure involved in
"resolution of unity". So the completeness of GCS has not yet been
established.

The main results are
\begin{itemize}
\item We show that our approach to construction of coherent states
works in the deformed case as well as in the classical one.
\item The well-known Arik -Coon oscillator
naturally arose in our approach as oscillator connected with the
$q$-Hermite polynomials in the same way as usual oscillator with standard
Hermite
polynomials.
\end{itemize}
\bigskip

\section{Coherent states for Rogers $\mathbf q$-Hermite polynomials
\hfill\break
$H_n(x;q)$ and Arik - Coon oscillator}
\subsection{$\mathbf q$-Oscillator, connected with  $H_n(x;q)$}
In the Hilbert space
\begin{equation}\label{RHS}
{\mathcal H}_q={\rm L}^2\left([-1,1];\,{\rm d}\mu_q(x)\right),
\qquad {\rm d}\mu_q(x)=\frac{(q;q)_{\infty}}{2\pi}\,
\frac{\left|(e^{2i\theta};q)_{\infty}\right|^2}{\sqrt{1-x^2}},
\end{equation}
we consider the orthonormal basis given by Rogers $q$-Hermite
polynomials
\begin{equation}\label{basis}
\varphi_n(x;q)={(q;q)_{n}^{-1/2}}H_n(x;q).
\end{equation}
The recurrence relations for the polynomials $\varphi_n(x;q)$ have
the form
\begin{equation}\label{rr}
x\varphi_n(x;q)=b_{n}\varphi_{n+1}(x;q)\!+b_{n-1}\varphi_{n-1}(x;q),\qquad
\varphi_0(x;q)=1,
\end{equation}
with coefficients
\begin{equation}\label{rrc}
b_{n}=\frac12\sqrt{1-q^{n+1}},\quad n\geq 0,\qquad b_{-1}=0.
\end{equation}
We consider ${\mathcal H}_q$ as the Fock space for the deformed
oscillator defined by the following operators
($\ket{n}:=\varphi_n(x;q)$)
\begin{gather}
X_q\ket{n}=b_{n}\ket{n+1}+b_{n-1}\ket{n-1};\\ P_q\ket{n}=i\left(
b_{n}\ket{n+1}-b_{n-1}\ket{n-1}\right);\\
a^{+}_q=(1-q)^{-\half}\left( X_q-iP_q \right),\quad
a^{+}_q\ket{n}=\sqrt{[n+1]_q} \ket{n+1};\\
a^{-}_q=(1-q)^{-\half}\left( X_q+iP_q \right),\quad
a^{-}_q\ket{n}=\sqrt{[n]_q} \ket{n-1},
\end{gather}
where
$$
[n]_q=\frac{1-q^n}{1-q}
$$
is the "mathematical" $q$-number and we consider the case
$0<q<1.$ From above relations we obtain
\begin{equation}\label{ACo}
a^{-}_qa^{+}_q-qa^{+}_qa^{-}_q=1,
\end{equation}
which means that we obtain the Arik - Coon oscillator.

The polynomials $\ket{n}=\varphi_n(x;q)$ are eigenfunctions of the
Hamiltonian
\begin{gather}
H_q={X_q}^2+{P_q}^2=a^{+}_qa^{-}_q+a^{-}_qa^{+}_q,\nonumber\\
H_q\ket{n}=\lambda_n\ket{n};\label{eve}\\
\lambda_0=\frac{4b_0^2}{1-q}=1;\quad
\lambda_n=\frac{4}{1-q}\left(b_{n-1}^2+b_{n}^2\right)=
\left([n]_q+[n+1]_q\right)\nonumber
\end{gather}
The equation (\ref{eve}) is equivalent to $q$-difference equation
for the Rogers $q$-Hermite polinomials
\begin{equation}
(1-q)D_q\left[w(x)D_q\varphi_n(x;q)\right]+
4q^{1-n}[n]_qw(x)\varphi_n(x;q)=0,
\end{equation}
where
\begin{equation}
w(x)={\left|(e^{2i\theta};q)_{\infty}\right|^2},\qquad
D_qf(x):=\frac{\delta_qf(x)}{\delta_qx},
\end{equation}
with
\begin{equation}
\delta_qf(e^{i\theta})=
f(q^{\half}e^{i\theta})-f(q^{-\half}e^{i\theta}),\quad
x=\cos\theta.
\end{equation}
The unitary equivalence of the position and the momentum operators
is given by a generalization of the Fourier transform. Besides,
the Hamiltonian (\ref{eve}) is invariant under generalized Fourier
transform as well as in the classical case.

Let's note that from the other point of view the connection of
Arik - Coon oscillator with Rogers $q$-Hermite polynomials was mentioned
in~\cite{Mac,FTV1}.

\subsection{Generalized Fourier transform for Rogers
$\mathbf q$-Hermite polynomials}

We define Generalized Fourier transform (GFT) related to the
orthonormal system $\left\{\varphi_n(x;q)\right\}_0^{\infty}$ by
the relation
\begin{equation}\label{}
{\rm F}_{\varphi}f(y)=
\int_ {-1}^{1} {{\rm K}_{\varphi}(x,y;-i)f(x)}{\mu_q(x)},
\end{equation}
where the Poisson kernel is given by
\begin{equation*}
{\rm K}_{\varphi}(x,y;t)\!=\!\Szi{n}t^n\varphi_n(x;q)\varphi_n(y;q)
\quad\text{and}\quad
{\rm K}_{\varphi}(x,y;-i)\!=
\!\lim_{t\rightarrow -i}{\rm K}_{\varphi}(x,y;t).
\end{equation*}
From q-Mehler formula for Rogers $q$-Hermite polynomials
$H_n(x;q)$ we obtain
\begin{equation}\label{}
{\rm K}_{\varphi}(x,y;t)=\frac{(t^2;q)_{\infty}}
{\left|(te^{2i\theta};q)_{\infty}(t;q)_{\infty} \right|^2}.
\end{equation}

It is easily to check, that our definition of GFT conforms with
given (from the different considerations) in~\cite{AAS}.

\subsection{Coherent states for the $\mathbf q$-oscillator,
connected with Rogers $\mathbf q$-Hermite polynomials $H_n(x;q)$}

The Barut - Girardello coherent states
are defined by
\begin{equation}\label{CS}
a^{-}_q\ket{z}=z\ket{z};\qquad
\ket{z}={\mathcal N}^{-1}\Szi{n}\frac{H_n(x;q)}{\sqrt{(q;q)_{n}}}
\,\,\frac{z^n}{\sqrt{[n]_q!}},
\end{equation}
where
\begin{equation}
{\mathcal N}^{2}=
\Szi{n}\frac{|z|^{2n}}{[n]_q!}={\tilde e}_q((1-q)|z|^{2});\qquad
R=\frac1{\sqrt{1-q}}.
\end{equation}
Here $R$ is the radius of convergence and by definition
\begin{equation}
{\tilde e}_q(x):=\Szi{n}\frac{x^n}{(q;q)_{n}}.
\end{equation}
(Let's note that we reserve the notation $e_q(x)$ for slightly
different definition of $q$-exponential function
\begin{equation}
e_q(x):=\Szi{n}\frac{x^n}{[n]_q!}={\tilde e}_q((1-q)x),
\end{equation}
which is more customary for quantum groups theory.)

Using the generation function for the Rogers $q$-Hermite
polinomials, from (\ref{CS}) we obtain
\begin{equation}
\ket{z}=\frac{{\tilde e}_q(e^{i\theta}\sqrt{1-q}z)
{\tilde e}_q(e^{-i\theta}\sqrt{1-q}z)}
{\sqrt{{\tilde e}_q((1-q)|z|^{2})}}.
\end{equation}
The overlap of two coherent states is given by
\begin{equation}
\bra{z_1}z_2\rangle_q={\tilde e}_q((1-q)\bar{z}_1z_2).
\end{equation}

To prove "resolution of identity" property
\begin{equation}
\label{RoI}
\int\int_{{\mathbb C}_{1/(1-q)}}\!\!\!\!\!\!\ket{z}
\raisebox{-3pt}{${}_q$}\langle{z}\vert
\text{d}\mu(|z|^2)\!=\!I,\quad
\text{d}\mu(|z|^2)\!=\!
{W}(|z|^2)\text{d}(\text{Re{z}})\text{d}(\text{Im {z}}),
\end{equation}
which means that the set of coherent states is complete (really,
over complete) we must solve the moment problem
\begin{equation}\label{}
\Izx{1/(1-q)}{t^n\frac{W(t)}{{\mathcal N}^{2}}}{t}=\frac{[n]_q!}{\pi}.
\end{equation}
Its solution is given by the distribution
\begin{equation}\label{}
W(t)=\frac{1-q}{\pi}\Szi{k}\frac{{\mathcal N}^{2}t}{e_q(qt)}
\delta\left(t-\frac{q^k}{1-q}\right),
\end{equation}
which can be obtained from the relation (\cite{AC}):
\begin{equation}\label{ji}
I_n(q)=\qIzx{1/(1-q)}{{e_q}^{-1}(qx)\,x^n}{x}=[n]_q!,
\end{equation}
where $I_n(q)$ is the so-called Jackson integral
\begin{equation}\label{ji1}
\qIzx{a}{f(x)}{x}:=a(1-q)\Szi{k}q^kf(q^ka).
\end{equation}
Thus the measure in (\ref{RoI}) is given by
\begin{align}
\text{d}\mu(|z|^2)&=
\frac{1-q}{\pi}\Szi{k}\frac{{\mathcal N}^{2}|z|^2}{e_q(q|z|^2)}
\delta\left(|z|^2-\frac{q^k}{1-q}\right)\text{d}^2z \nonumber\\
{}&=\frac{1-q}{\pi}\Szi{k}\frac{|z|^2e_q(|z|^2)}{e_q(q|z|^2)}
\delta\left(|z|^2-\frac{q^k}{1-q}\right)\text{d}^2z. \label{qH41}
\end{align}

For the completeness we give the simple proof of the essential
relation (\ref{ji}). It is well known that for $q$-derivative
defined as
\begin{equation}\label{p1}
\qd f(x)=\frac{f(x)-f(qx)}{x(1-q)}\,,
\end{equation}
one has (and easily verified)
\begin{equation}\label{p2}
\qd x^n=[n]_qx^{n-1}\,,\qquad\quad \qd e_q(x)=e_q(x)\,.
\end{equation}
From the Leibnitz relation for $q$-derivative
\begin{equation}\label{p3}
\qd\left(u(x)v(x)\right)=u(x)\qd v(x)+v(qx)\qd u(x)\,,
\end{equation}
it follows the integration by parts rule for the Jackson integral
in the form
\begin{equation}\label{p4}
\qIzx{a}{\!u(x)\qd \!v(x)}{x}\!=\!\qd\!\left(u(x)v(x)\right)
\rule[-7pt]{0.5pt}{25pt}^{\,a}_{\,0}
-\qIzx{a}{\!v(qx)\qd \!u(x)}{x}.
\end{equation}
By differentiating the identity  $e_q(x)\cdot e_q^{\,-1}(x)=1$ one
obtains
\begin{equation}\label{p5}
\qd e_q^{\,-1}(x)=-e_q^{\,-1}(qx)\,.
\end{equation}
Using this relation and integrating by parts the Jackson integral $I_n(q)$
(\ref{ji}) we have
\begin{align}
I_n(q)&=\qIzx{\frac{1}{1-q}}{x^n\left(-\qd e_q^{\,-1}(x)\right)}{x}=
\nonumber\\
{}&=x^n\left(-\qd e_q^{\,-1}(x)\right)
\rule[-7pt]{0.5pt}{25pt}^{\frac{1}{1-q}}_{\,0}+
\qIzx{\frac{1}{1-q}}{e_q^{\,-1}(qx)\qd x^n}{x}\,.\label{p6}
\end{align}
In view of the relation $e_q^{\,-1}(\frac{1}{1-q})=0$ one gets from
(\ref{p6}) the recurrent relation
\begin{equation}\label{p7}
I_n(q)=[n]_q\,I_{n-1}(q)\,, \quad n\geq 1\,,
\end{equation}
so that
\begin{equation}\label{p8}
I_n(q)=[n]_q!\,I_{0}(q)\,.
\end{equation}
which gives the  relation (\ref{ji}), because
\begin{equation}\label{p9}
I_{0}(q)=\qIzx{\frac{1}{1-q}}{e_q^{\,-1}(qx)}{x}=
- e_q^{\,-1}(x)\,\rule[-7pt]{0.5pt}{25pt}^{\frac{1}{1-q}}_{\,0}=
e_q^{\,-1}(0)=1.
\end{equation}

Let us note that over completeness of this coherent states was
considered from another point of view, for example in ~\cite{Per}.

\section{Coherent states for discrete $\mathbf q$-Hermite polynomials II}
\subsection{Discrete $\mathbf q$-Hermite polynomials II}
The discrete $q$-Hermite polynomials of II-type are defined by
\begin{equation}\label{qHIIa}
{\tilde h}_n(x;q)=x^n \sfito{q^{-n}}{q^{-n+1}}{0}{q^2;\,-x^{-2}},
\end{equation}
where
\begin{equation}
\sfito{a}{b}{c}{q;z}=
\Szi{k}\frac{(a;q)_k(b;q)_k}{(c;q)_k}\,\frac{z^k}{(q;q)_k}.
\end{equation}
These polynomials fulfill the following orthogonality condition
\begin{equation}\label{oc}
c(1\!-\!q)\!\!\ \sum_{k=-\infty}^{\infty}
 \!\!\!\left[\! {\tilde h}_m(cq^k;q){\tilde h}_n(cq^k;q)\!+\!
{\tilde h}_m(-cq^k;q){\tilde h}_n(-cq^k;q)\!
\right]\!w(cq^k)q^k\!=\!0,
\end{equation}
where $m\neq n$, $w(x)=[(ix;q)_{\infty}(-ix;q)_{\infty}]^{-1}$ and
$c>0$.
We denote by ${\mathcal H}$ the Hilbert space, spanned by the set
of discrete $q$-Hermite polynomials of II-type. The orthonormal
polynomials
\begin{equation}
\psi_n(x;q)={\tilde h}_n(x;q)q^{\half n^2}{(q;q)_n}^{-\half},
\end{equation}
form the basis in ${\mathcal H}$ and fulfill the recurrence
relation
\begin{equation}
x\psi_n(x;q)= b_{n}\psi_{n+1}(x;q)+b_{n-1}\psi_{n-1}(x;q);\qquad
\psi_0(x;q)=1,
\end{equation}
where
\begin{equation}\label{}
b_{n}=q^{-n-\half}\sqrt{1-q^{n+1}},\quad n\geq 0;\qquad b_{-1}=0.
\end{equation}

\subsection{Deformed oscillator connected with discrete
$\mathbf q$-Hermite polynomials II}

We consider ${\mathcal H}_q$ as the Fock space for a deformed
oscillator defined by the following operators
($\ket{n}:=\psi_n(x;q)$)
\begin{gather}
X_q\ket{n}=b_{n}\ket{n+1}+b_{n-1}\ket{n-1};\\ P_q\ket{n}=i\left(
b_{n}\ket{n+1}-b_{n-1}\ket{n-1}\right);\\
a^{+}_q=\half\sqrt{\frac{q}{1-q}}\left( X_q-iP_q \right),\quad
a^{+}_q\ket{n}=\sqrt{\frac{q}{1-q}}b_{n}\ket{n+1};\\
a^{-}_q=\half\sqrt{\frac{q}{1-q}}\left( X_q+iP_q \right),\quad
a^{-}_q\ket{n}=\sqrt{\frac{q}{1-q}}b_{n-1}\ket{n-1},\\
N\ket{n}=n\ket{n};\quad a^{-}_qa^{+}_q=q^{-2N}[N+I]_q; \quad
a^{+}_qa^{-}_q=q^{-2N+2I}[N]_q;
\end{gather}
where $[n]_q=\frac{1-q^n}{1-q}$ is the "mathematical" $q$-number.
From the above relations we obtain
\begin{equation}\label{ACo2}
a^{-}_qa^{+}_q-q^{-1}a^{+}_qa^{-}_q=q^{-2N},\quad\text{or}\quad
a^{-}_qa^{+}_q-q^{-2}a^{+}_qa^{-}_q=q^{-N}.
\end{equation}
The polynomials $\ket{n}=\varphi_n(x;q)$ are eigenfunctions of the
Hamiltonian
\begin{gather}
H_q=\half\frac{q}{1-q}\left({X_q}^2+{P_q}^2\right)
=a^{+}_qa^{-}_q+a^{-}_qa^{+}_q,\\
H_q\ket{n}=\lambda_n\ket{n};\label{eve2}\\
\lambda_n=
q^{-2n}[n+1]_{q}+q^{2-2n}[n]_{q},\quad n\geq 0.
\end{gather}
The equation (\ref{eve2}) is equivalent to $q$-difference equation
for the discrete $q$-Hermite polynomials II
\begin{gather}
-(1-q^n)x^2{\tilde h}_n(x;q)=\\
 =q{\tilde h}_n(x-i;q)-
(1+q+x^2){\tilde h}_n(x;q)+(1+x^2){\tilde h}_n(x+i;q).
\end{gather}

\subsection{Coherent states for the $\mathbf q$-oscillator,
connected with  ${\tilde h}_n(x;q)$}

The Barut - Girardello coherent states are defined by
\begin{equation}\label{CS2}
a^{-}_q\ket{z}\!=\!z\ket{z};\quad \ket{z}\!=\!{\mathcal
N}^{-1}\Szi{n} \left(q(1-q)\right)^{\half n}q^{n^2-n}{\tilde
h}_n(x;q) \,\frac{z^n}{(q;q)_n},
\end{equation}
where
\begin{equation}
{\mathcal N}^{2}=\Szi{n}\left(\frac{1-q}{q}\right)^{n}q^{n^2}
\frac{z^{2n}}{(q;q)_n}\equiv\epsilon_b(z^{2});\qquad R=\infty.
\end{equation}
Using the generation function for the for the discrete $q$-Hermite
polynomials II, from (~\cite{KS}) we obtain

\begin{equation}
\ket{z}=\frac{1}{\epsilon_b(z^{2})} \left({i\sqrt{q(1-q)}z;q}
\right)_{\infty}
\sfioo{i{x}}{i\sqrt{q(1-q)}z}{q;-i\sqrt{q(1-q)}z},
\end{equation}

where

\begin{equation}
\sfioo{a}{b}{q;z}
=\Szi{k}\frac{(a;q)_k}{(b;q)_k}\,\,\frac{z^k}{(q;q)_k}.
\end{equation}

\subsection*{Acknowledgements}
This research was supported by RFFI grants No 03-01-00837, 03-01-00593

\end{document}